\documentclass[12pt,letter]{article}

\usepackage{mystyle}
\usepackage{mythmstyle_report}

\begin{document}
\title{Structure of the automorphism group of the augmented cube graph}

\author{Ashwin~Ganesan%
  \thanks{53 Deonar House, Deonar Village Road, Deonar, Mumbai 400088, Maharashtra, India. Email:  
\texttt{ashwin.ganesan@gmail.com}.}
% A.G.~is the corresponding author.}
}

\date{}
\vspace{10cm}

\maketitle

\begin{abstract}
\noindent  The augmented cube graph $AQ_n$ is the Cayley graph of $\mathbb{Z}_2^n$ with respect to the set of $2n-1$ generators $\{e_1,e_2, \ldots,e_n, 00\ldots0011, 00\ldots0111, 11\ldots1111 \}$. It is known that the order of the automorphism group of the graph $AQ_n$ is $2^{n+3}$, for all $n \ge 4$.  In the present paper, we obtain the structure of the automorphism group of $AQ_n$ to be 
\[ \Aut(AQ_n) \cong \mathbb{Z}_2^n \rtimes D_8~~(n \ge 4),\] 
where $D_8$ is the dihedral group of order 8.  It is shown that the Cayley graph $AQ_3$ is non-normal and that $AQ_n$ is normal for all $n \ge 4$.  We also  analyze the clique structure of $AQ_4$ and show that the automorphism group of $AQ_4$ is isomorphic to that of $AQ_3$:  \[ \Aut(AQ_4) \cong \Aut(AQ_3) \cong (D_8 \times D_8) \rtimes C_2.\] 
All the nontrivial blocks of $AQ_4$ are also determined. 
\end{abstract}

\bigskip
\noindent\textbf{Index terms} --- augmented cubes; automorphisms of graphs; normal Cayley graphs; clique structure; nontrivial block systems.

%\newpage
%---------------------------------------------------------------

\vspace{+0.5cm}
%---------------------------------------------------------------
\section{Introduction}

Cayley graphs have been well-studied as a topology for interconnection networks due to their high symmetry, low diameter, high connectivity, and embeddable properties.  Perhaps the most well-studied topology is the hypercube graph, for which there are many equivalent definitions.  One definition is to view the hypercube as a Cayley graph.

Given a group $H$ and a subset $S \subseteq H$, the Cayley graph of $H$ with respect to $S$, denoted $\Cay(H,S)$, is the graph whose vertex set is $H$, and with arc set $\{ (h, sh): h \in H, s \in S \}$.  If the identity group element $e$ is not in $H$, then $\Cay(H,S)$ has no self-loops, and if $H$ is closed under inverses, then $(h,g)$ is an arc iff $(g,h)$ is an arc; in this  case, the Cayley graph can be viewed as an undirected graph.  The Cayley graph $\Cay(H,S)$ is connected iff $S$ generates $H$.

Let $\mathbb{Z}_2^n$ denote the direct product of $n$ copies of the cyclic group of order 2.  Let $e_i$ denote the unit vector in the vector space $\mathbb{Z}_2^n$ that has a 1 in the $i$th coordinate and zero in the remaining coordinates. The hypercube graph is defined to be the Cayley graph of $\mathbb{Z}_2^n$ with respect to the set of generators $S = \{e_1,e_2, \ldots, e_n\}$.  

The hypercube graph has been well-studied in the literature as a topology for interconnection networks.  In order to improve its connectivity and fault-tolerance properties, the folded hypercube graph was proposed \cite{ElAmawy:Latifi:1991}.   The folded hypercube graph is obtained by taking the hypercube graph and joining each vertex $x$ to the unique vertex $\overline{x}$ that is farthest away from $x$ in the hypercube.  Thus, the folded hypercube is the Cayley graph of $\mathbb{Z}_2^n$ with respect to the set of generators $S = \{e_1,e_2,\ldots,e_n, 1\cdots11\}$.  The addition of these edges to the hypercube graph yields a new graph whose diameter is smaller.

The augmented cube graph was introduced in Choudum and Sunitha \cite{Choudum:Sunitha:2002} as a topology for consideration in interconnection networks.  The augmented cube graph of dimension $n$, denoted $AQ_n$, is the Cayley graph of $\mathbb{Z}_2^n$ with respect to the set of $2n-1$ generators $\{e_1,e_2,\ldots,e_n\}$ $ \dot{\cup} \{00\cdots0011, $ $00\cdots0111, 00\cdots1111, 11\cdots1111\}$. Recall that the hypercube graph is the Cayley graph of $\mathbb{Z}_2^n$ with respect to the generator set $\{e_1,\ldots,e_n\}$.  The additional $n-1$ generators in the definition of $AQ_n$ augment the edge set of the hypercube graph, so that the augmented graph has smaller diameter and better connectivity and embeddable properties. For equivalent definitions of the augmented cube graph and a study of some of its properties, see  \cite{Choudum:Sunitha:2002} \cite{Choudum:Sunitha:2008} \cite{Choudum:Sunitha:2001}.

Let $X=(V,E)$ be a simple, undirected graph.  The automorphism group of $X$ is defined to be the set of all permutations of the vertex set that preserves adjacency: $\Aut(X):=\{g \in \Sym(V): E^g = E\}$.    While one can often obtain some automorphisms of a graph, it is often difficult to prove that one has obtained the (full) automorphism group. Automorphism groups of graphs that have been proposed as the topology of interconnection networks have been investigated by several authors; see \cite{Deng:Zhang:2011} \cite{Deng:Zhang:2012} \cite{Ganesan:JACO} \cite{Ganesan:DM:2013}  \cite{Feng:2006} \cite{Mirafzal:2011} \cite{Zhang:Huang:2005} \cite{Zhou:2011}.

Let $G:=\Aut(AQ_n)$ denote the automorphism group of the augmented cube graph.  Since $AQ_n = \Cay(\mathbb{Z}_2^n,S)$ is a Cayley graph, its automorphism group can be expressed as a rotary product $G = \mathbb{Z}_2^n \times_{rot} G_e$, where $G_e$ is the stabilizer in $G$ of the identity vertex $e$ (cf. Jajcay \cite{Jajcay:1994} \cite{Jajcay:2000}).  In  \cite{Choudum:Sunitha:2008} \cite{Choudum:Sunitha:2001}, it was shown that the order of the automorphism group of $AQ_n$ is exactly $2^{n+3}$, for all $n \ge 4$.  

In the present paper, the structure of the automorphism group $G:=\Aut(AQ_n)$ is determined to be $G \cong \mathbb{Z}_2^n \rtimes D_8$, for all $n \ge 4$.  Hence, in addition to proving that $G_e \cong D_8$, we also strengthen the rotary product to a semidirect product.  The clique structure of $AQ_4$ is investigated further, and we show that the automorphism group of $AQ_4$ is isomorphic to that of $AQ_3$: $\Aut(AQ_4) \cong \Aut(AQ_3) \cong (D_8 \times D_8) \rtimes C_2$.  We also determine all the nontrivial blocks of $AQ_4$.

An open problem in the literature is to determine, given a group $H$ and a subset $S \subseteq H$, whether $\Cay(H,S)$ is normal \cite{Xu:1998}.  Every Cayley graph $\Cay(H,S)$ admits $R(H) \Aut(H,S)$ as a subgroup of automorphisms; here  $R(H)$ denotes the right regular representation of $H$,  and $\Aut(H,S)$ is the set of automorphisms of $H$ that fixes $S$ setwise (cf. \cite{Biggs:1993} \cite{Godsil:Royle:2001}).  A Cayley graph $\Cay(H,S)$ is said to be normal if its full automorphism group is $R(H) \Aut(H,S)$.  In the present paper, we prove that the Cayley graph $AQ_3$ is non-normal, and we prove that $AQ_n$ is normal for all $n \ge 4$.

\bigskip \noindent \textbf{Notation.} We mention the notation used in this paper. Given a graph $X=(V,E)$ and a vertex $v \in V$, $X_i(v)$ denotes the set of vertices of $X$ whose distance to vertex $v$ is exactly $i$. Thus, if $X=\Cay(H,S)$, then $X_1(e)=S$. The subgraph of $X$ induced by a subset of vertices $W \subseteq V$ is denoted $X[W]$, and so $X[X_i(v)]$ denotes the subgraph of $X$ induced by the $i$th layer of the distance partition of $X$ with respect to vertex $v$.   The dihedral group of order $2n$ is denoted $D_{2n}$. 

Given a Cayley graph $X = \Cay(H,S)$, $e$ denotes the identity element of the group $H$ and also denotes the corresponding vertex of $X$.  Throughout this paper, $H$ will be the abelian group $\mathbb{Z}_2^n$.  In this case, $R(H) = R(\mathbb{Z}_2^n)$ is the set $\{\rho_z: z \in \mathbb{Z}_2^n\}$ of translations by $z$, i.e. $\rho_z: \mathbb{Z}_2^n \rightarrow \mathbb{Z}_2^n, x \mapsto x+z$, and $\rho_z$ is an automorphism of the Cayley graph.  We also let $\mathbb{Z}_2^n$ denote the vector space $\mathbb{F}_2^n$. 

Let $G:=\Aut(X)$.  Then $G_e$ is the set of automorphisms of $X$ that fixes the vertex $e$, and $L_e$ denotes the set of automorphisms of $X$ that fixes the vertex $e$ and each of its neighbors. 

%-----------------------------------------------------------------

\section{$AQ_2$ and $AQ_3$}

\begin{Lemma}
 Let $AQ_2$ be the augmented cube graph $\Cay(\mathbb{Z}_2^2,S)$, where $S = \{e_1,e_2,11\}$. Then $\Aut(AQ_2) \cong S_4$ and $AQ_2$ is normal.
\end{Lemma}

\noindent \emph{Proof}:
Note that $AQ_2$ is the complete graph on 4 vertices. Hence, $\Aut(AQ_2) \cong S_4$ and the stabilizer in $G:=\Aut(AQ_2)$ of the identity vertex $e$ is $G_e \cong S_3$.  Recall that the Cayley graph $\Cay(\mathbb{Z}_2^2,S)$ is normal iff $G_e \subseteq \Aut(\mathbb{Z}_2^n,S)$ (cf. Xu \cite[Proposition 1.5]{Xu:1998}). The set $\Aut(\mathbb{Z}_2^3,S)$ of automorphisms of the group $\mathbb{Z}_2^2$ that fixes $S$ setwise is exactly the set of invertible linear transformations of the vector space $\mathbb{Z}_2^2$ that fixes $S$ setwise. The elements of $\Aut(\mathbb{Z}_2^2,S)$ are uniquely determined by their action on the basis $\{e_1,e_2\}$ of the vector space. Let $\phi \in \Aut(\mathbb{Z}_2^2,S)$. The image of $\phi(e_1)$ can be chosen from $S$ in 3 ways, and the image $\phi(e_2)$ can be chosen from $S - \{\phi(e_1)\}$ in 2 ways. Hence, $|\Aut(\mathbb{Z}_2^2,S)|=6$, which equals $|G_e|$.  Thus, $AQ_2$ is normal, i.e. it has the smallest possible full automorphism group $R(\mathbb{Z}_2^2) \Aut(\mathbb{Z}_2^2,S) \cong \mathbb{Z}_2^2 \rtimes S_3$.
\qed

\begin{Theorem}
 Let $AQ_3$ be the Cayley graph $\Cay(\mathbb{Z}_2^3,S)$, where $S:=\{e_1,e_2,e_3,011,111\}$. Let $G:=\Aut(AQ_3)$. Then $G_e \cong D_8 \times C_2$, $\Aut(\mathbb{Z}_2^3,S) \cong D_8$, and $AQ_3$ is non-normal.
\end{Theorem}

\noindent \emph{Proof}:
Let $X:=AQ_3$ and $G:=\Aut(AQ_3)$.  We show that the stabilizer in $G$ of the vertex $e$ is $G_e \cong D_8 \times C_2$ and that $\Aut(\mathbb{Z}_2^3,S) \cong D_8$, whence $\Aut(\mathbb{Z}_2^3,S) \ne G_e$ and $AQ_3$ is non-normal.

\begin{figure}
\centering
\begin{tikzpicture}[scale=1.5]
\vertex[fill] (e) at (0,0) [label=below:$e$] {};

\vertex[fill] (v100) at (2,-1) [label=below:$111$] {};
\vertex[fill] (v111) at (2,1) [label=above:$e_1$] {};

\vertex[fill] (v011) at (3,0) [label=right:$011$] {};

\vertex[fill] (v010) at (4,-1) [label=below:$e_2$] {};
\vertex[fill] (v001) at (4,1) [label=above:$e_3$] {};

\vertex[fill] (v110) at (6,-1) [label=below:$101$] {};
\vertex[fill] (v101) at (6,1) [label=above:$110$] {};

\node at (3,-2) {${X_1(e)}$};
\node at (6,-2) {${X_2(e)}$};
\path
    (v100) edge (v011)
    (v111) edge (v011)
    (v010) edge (v011)
    (v001) edge (v011)

    (v100) edge (v111)
    (v001) edge (v010)

    (v110) edge (v101)
;
\end{tikzpicture}
\caption{Distance partition of $AQ_3$ (edges between different layers are not drawn).} 
\label{fig:distance:partition:AQ3}
\end{figure}
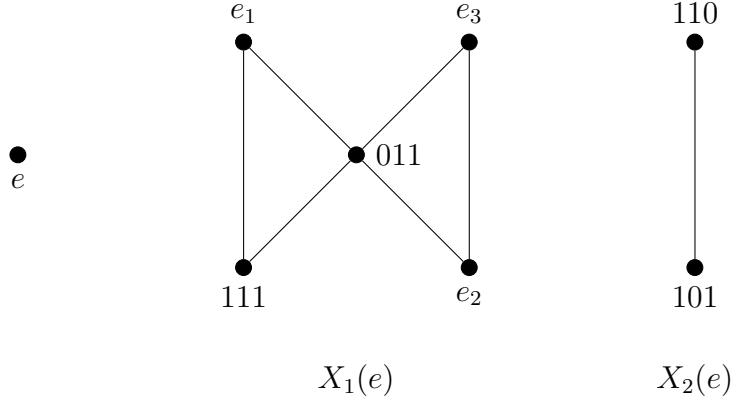

The distance partition of $AQ_3$ with respect to the identity vertex $e$ is shown in Figure~\ref{fig:distance:partition:AQ3}.  The vertex $e$ is adjacent to each vertex in $X_1(e)$ and each $u \in X_1(e) - \{011\}$ is adjacent to both vertices in $X_2(e)$ (these edges between the layers are not shown in the figure). The graph $AQ_3$ is 5-regular. 

The subgraph of $X$ induced by the $i$th layer of the distance partition is denoted $X[X_i(e)]$ .  The stabilizer $G_e$ acts on the subgraph $X[X_i(e)]$ as a group of automorphisms. Observe from Figure~\ref{fig:distance:partition:AQ3} that an automorphism $g \in G_e$ can act on $X_1(e)$ by permuting the endpoints of the edge $\{e_1,111\}$, by independently permuting the endpoints of the edge $\{e_2,e_3\}$, and by thereafter possibly also interchanging these two edges.  Hence, the automorphism group of the induced subgraph $X[X_1(e)]$ is that of two independent edges (or its complement graph - the square) and hence is isomorphic to $D_8$. 

The automorphism group of the induced subgraph $X[X_2(e)]$ is clearly $C_2$. Since each vertex in $X_2(e) - \{011\}$ is adjacent to all vertices in $X_1(e)$ and since $G_e$ fixes the vertex $011$, $G_e  \cong D_8 \times C_2$.  This also proves that $\Aut(AQ_3) = R(\mathbb{Z}_2^3) \times_{rot} (D_8 \times C_2)$. 

We now determine $\Aut(\mathbb{Z}_2^3,S)$. Recall that the automorphisms of the group $\mathbb{Z}_2^n$ are precisely the invertible linear transformations of the vector space $\mathbb{Z}_2^n$ (cf. \cite[p. 136]{Dummit:Foote:2004}).  Hence $\Aut(\mathbb{Z}_2^3,S)$ is precisely the set of elements of $G_e$ that can be realized as linear transformations of the vector space $\mathbb{Z}_2^3$.  We first show that $\Aut(X[X_1(e)]) \cong \Aut(\mathbb{Z}_2^3,S)$.  
Each automorphism $\phi$ of the induced subgraph $X[X_1(e)]$ also specifies the image $\phi(e_i)$ of the vertex $e_i~~ (i=1,2,3)$. The resulting image of the basis $\{e_1,e_2,e_3\}$ extends uniquely to a linear transformation $\phi'$ of the vector space $\mathbb{Z}_2^3$. Let $f: \Aut(X[X_1(e)]) \rightarrow \Aut(\mathbb{Z}_2^3,X_1(e)), \phi \mapsto \phi'$ be the map that takes an automorphism of the induced subgraph $X[X_1(e)]$ to its linear extension $\phi'$.  We show that the map $f$ is well-defined and a bijective homomorphism.

Let $\phi \in \Aut(X[X_1(e)])$. Observe from Figure~\ref{fig:distance:partition:AQ3} that $\phi$ fixes the vertex $011$ and permutes the vertices $X_1(e)-\{011\} = \{e_1,e_2,e_3,111\}$ among themselves.  In particular, if $\phi$ takes the basis element $e_i$ to $\phi(e_i)$ $(i=1,2,3)$, then $\phi$ must take the fourth vertex $111=e_1+e_2+e_3$ to the remaining vertex (i.e. to the unique vertex in the singleton set $\{e_1,e_2,e_3,111\} - \{\phi(e_1), \phi(e_2), \phi(e_3) \}$. This remaining vertex is the sum $\phi(e_1)+\phi(e_2)+\phi(e_3)$ because every 3 vectors in the set $\{e_1,e_2,e_3,111\}$ sum to the fourth vector in the set.  Also, the subsets $\{e_2,e_3\}$ and $\{e_1,111\}$ are blocks of $\phi$, and the sum of the two vertices in each subset is $011$.  In other words, if $\phi \in \Aut(X[X_1(e)])$, then $\phi$ acts linearly on the vertices in $X_1(e)$, and so the extension of $\phi$ to a linear transformation $\phi'$ is well-defined.  

An automorphism of the induced subgraph $X([X_1(e)])$ must take vertices $e_1,e_2$ and $e_3$ to 3 vertices in $\{e_1,e_2,e_3,111\}$. Any 3 of these 4 vectors are linearly independent, whence the linear extension of $\phi$ is an invertible linear transformation. Thus, the map $f$ is well-defined.

Different automorphisms of $X([X_1(e)])$ induce different images of the basis $\{e_1,e_2,e_3\}$ and hence different linear extensions.  Thus, $f$ is injective. Since $\Aut(\mathbb{Z}_2^3,S) \subseteq G_e$, every element $\phi' \in \Aut(\mathbb{Z}_2^3,S)$, when restricted to $S$, is an automorphism $\phi$ of $X([X_1(e)])$.  The linear extension of $\phi$ is unique and hence equals $\phi'$. Thus, $f$ is surjective. The linear extension map is also a homomorphism.  Hence $f$ is an isomorphism.

Since $\Aut(X[X_1(e)]) \cong D_8$, we have that $\Aut(\mathbb{Z}_2^3,S) \cong D_8$. Thus, $\Aut(\mathbb{Z}_2^3,S)$ is a proper subgroup of $G_e = D_8 \times C_2$ and $AQ_3$ is non-normal.
\qed

\begin{Lemma}
 Let $AQ_3$ denote the augmented cube graph $\Cay(\mathbb{Z}_2^3,S)$ of dimension 3. Then $\Aut(AQ_3) \cong (D_8 \times D_8) \rtimes C_2$. 
\end{Lemma}

\noindent \emph{Proof}:
Since $|S|=5$, $AQ_3$ is a 5-regular graph on 8 vertices. Its complement graph $\overline{AQ_3}$ is $\Cay(\mathbb{Z}_2^3,S')$, where $S':=\mathbb{Z}_2^3-S-\{0\} = \{110,101\}$. It is easy to see that $\overline{AQ_3}$ is the union of two disjoint squares. Hence $\Aut(\overline{AQ_3}) \cong (D_8 \times D_8) \rtimes C_2$.
\qed

%------------------------------------------

\section{Automorphism group and normality of $AQ_n$ $(n \ge 4)$}

Let $AQ_n$ be the augmented cube graph $\Cay(\mathbb{Z}_2^n,S)$ $(n \ge 4)$. Let $G:=\Aut(AQ_n)$. Since $AQ_n$ is a Cayley graph, its automorphism group is the rotary product $R(\mathbb{Z}_2^n) \times_{rot} G_e$. It was shown in \cite{Choudum:Sunitha:2001}~\cite{Choudum:Sunitha:2008} that the number of automorphisms of the graph $AQ_n$ that fix the identity vertex $e$ is exactly 8, i.e. $|G_e|=8$. Thus, it was shown that the order of the automorphism group of $AQ_n$ is $2^n \times 8 = 2^{n+3}$.  In the present section, we determine the structure of the automorphism group of $AQ_n$: we show that $G_e \cong D_8$ and that every one of these 8 automorphisms can be realized as a linear transformation of the vector space $\mathbb{Z}_2^n$. It then follows that $G_e \subseteq \Aut(\mathbb{Z}_2^n,S)$ and that the rotary product expression can be strengthened to the semidirect product $G = R(\mathbb{Z}_2^n) \rtimes D_8$. 

We first recall the following result from the literature:

\begin{Lemma} \cite{Choudum:Sunitha:2001} \cite{Choudum:Sunitha:2008}
 Let $X:=AQ_n = \Cay(\mathbb{Z}_2^n,S)$ be the augmented cube graph $(n \ge 4)$ and $G:=\Aut(X)$. Then $|G_e|=8$.
\end{Lemma}

We recall briefly the proof given in \cite{Choudum:Sunitha:2001} for the result $|G_e|=8$.  First, it can be shown that the induced subgraph $X[X_1(e)]$ has exactly 8 automorphisms (cf. \cite[Lemma 3.2]{Choudum:Sunitha:2001}).  It can also be shown that every automorphism of the graph $AQ_n$ that fixes the identity vertex $e$ and each of its neighbors is trivial, i.e. $L_e=1$ (cf. \cite[Theorem 3.1]{Choudum:Sunitha:2001}).  It follows that $|G_e| \le 8$.  Finally, Table 1 in \cite{Choudum:Sunitha:2001} defines a set of 8 permutations of the vertex set of $AQ_n$ that fix the identity vertex $e$ and that are automorphisms of $AQ_n$.  Hence $|G_e|=8$.  The paper \cite{Choudum:Sunitha:2001} doesn't show the verification that the eight permutations given in \cite[Table 1]{Choudum:Sunitha:2001} preserve adjacency of $AQ_n$.  In the present section, we express these 8 permutations as linear extensions of automorphisms of the induced subgraph $X[X_1(e)]$.

\begin{Proposition} \label{prop:normal:AQn}
Let $X = AQ_n (n \ge 4)$ and $G:=\Aut(X)$. Then each of the 8 elements in $G_e$ can be realized as a linear transformation  of the vector space $\mathbb{Z}_2^n$ that fixes $S$ setwise, and so $G_e \subseteq \Aut(\mathbb{Z}_2^n,S)$. 
\end{Proposition}

\begin{figure}
\centering
\begin{tikzpicture}[scale=2]

\vertex[fill] (ve1) at (0,0) [label=below:$e_1$] {};
\vertex[fill] (ve2) at (1,0) [label=below:$e_2$] {};
\vertex[fill] (ve3) at (2,0) [label=below:$e_3$] {};

\vertex[fill] (ven_2) at (5,0) [label=below:$e_{n-2}$] {};
\vertex[fill] (ven_1) at (6,0) [label=below:$e_{n-1}$] {};

\vertex[fill] (v11dot1) at (-0.5,1) [label=above:$11\cdot\cdot1$] {};
\vertex[fill] (v01dot1) at (0.5,1) [label=above:$01\cdot\cdot1$] {};
\vertex[fill] (v001dot1) at (1.5,1) [label=above:$001\cdot\cdot1$] {};
\vertex[fill] (v0001dot1) at (2.5,1) [label=above:$0001\cdot\cdot1$] {};

\vertex[fill] (v0dot111) at (4.5,1) [label=above:$0\cdot\cdot111$] {};
\vertex[fill] (v0dot011) at (5.5,1) [label=above:$00\cdot\cdot011$] {};
\vertex[fill] (ve_n) at (6.5,1) [label=above:$e_n$] {};

\node at (3.5,0.5) {$\cdots$};

\path
    (ve1) edge (v11dot1)
    (ve1) edge (v01dot1)

    (ve2) edge (v01dot1)
    (ve2) edge (v001dot1)

    (ve3) edge (v001dot1)
    (ve3) edge (v0001dot1)

    (ven_2) edge (v0dot111)
    (ven_2) edge (v0dot011)

    (ven_1) edge (v0dot011)
    (ven_1) edge (ve_n)

    (v11dot1) edge (v01dot1)
    (v01dot1) edge (v001dot1)
    (v001dot1) edge (v0001dot1)
    (v0dot111) edge (v0dot011)
    (v0dot011) edge (ve_n)

;
\end{tikzpicture}
\caption{Induced subgraph of the first layer of the distance partition of $AQ_n$ $(n \ge 4)$.} 
\label{fig:induced:subgraph:X1e:AQn}
\end{figure}
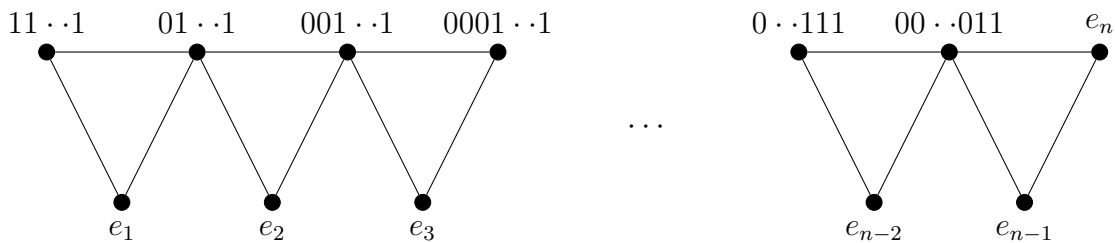

\noindent \emph{Proof}:
The induced subgraph $X[X_1(e)]$ is shown in Figure~\ref{fig:induced:subgraph:X1e:AQn}, where the subgraph is drawn horizontally. It can be seen from Figure~\ref{fig:induced:subgraph:X1e:AQn} that automorphisms of the induced subgraph $X[X_1(e)]$ are realized by permuting the endpoints of the edge $\{e_1,11\cdots1\}$, by independently permuting the endpoints of the edge $\{e_{n-1},e_n\}$, and by possibly also interchanging these two edges (in which case the automorphism flips the graph left to right).  Hence, $\Aut(X[X_1(e)]) \cong D_8$.

Define $f_1$ to be the map that takes vertex $e_1$ to $11\cdots1$, and $e_i$ to $e_i (i=2,\ldots,n)$. Let $f_1'$ be the linear extension of $f_1$.  Then $f_1'$ takes $11\cdots1=e_1+\cdots+e_n$ to $f_1'(e_1)+\cdots+f_1'(e_n)=11\cdots1+e_2+\cdots+e_n=e_1$.  Similarly, $f_1'$ takes $e_k+e_{k+1}+\cdots+e_n$ $(k \ge 2)$ to itself.  Thus, $f_1'$ just interchanges the endpoints of the edge $\{e_1,11\cdots1\}$ in the induced subgraph $X[X_1(e)]$.  Since $f_1'$ is an invertible linear transformation and fixes $X_1(e)$ setwise, $f_1' \in \Aut(\mathbb{Z}_2^n,S) \subseteq G_e$; hence, $f_1'$ is an automorphism of the graph $X$.  We have thus realized an element $f_1'$ of $G_e$ as a linear transformation of the vector space $\mathbb{Z}_2^n$.

Similarly, define $f_2$ to be the map that takes $e_{n-1}$ to $e_n$, $e_n$ to $e_{n-1}$, and $e_i$ to $e_i$ $(i=1,2,\ldots,n-2)$. Define $f_3$ to be the map that interchanges $e_1$ and $e_{n-1}$, interchanges $e_2$ and $e_{n-2}$, and so on, and takes $e_n$ to $11\cdots1$. Let $f_2'$ and $f_3'$ be the linear extensions of $f_2$ and $f_3$, respectively. It can be seen that the group of linear transformations $\langle f_1',f_2',f_3' \rangle$ is contained in $G_e$ and is isomorphic to $\langle f_1,f_2,f_3 \rangle \cong D_8$.  Thus, $G_e \cong D_8$ and $G_e \subseteq \Aut(\mathbb{Z}_2^n,S)$.
\qed

We have proved above that

\begin{Theorem} \label{thm:autgroup:sdp:AQn}
 Let $X=AQ_n = \Cay(\mathbb{Z}_2^n,S) (n \ge 4)$ and $G:=\Aut(X)$. Then $G_e \cong D_8$, $AQ_n$ is a normal Cayley graph, and $\Aut(AQ_n) \cong \mathbb{Z}_2^n \rtimes D_8$. 
\end{Theorem}

%-------------------------------------------------------------
\section{Automorphism group of $AQ_4$}

In this section, we study the clique structure of $X:=AQ_4 = \Cay(\mathbb{Z}_2^4,S)$, where $S:=\{e_1,e_2,e_3,e_4,0011,0111,1111\}$. In particular, we describe the action of the automorphism group of $AQ_4$ on the set of all its maximum cliques.  We show that the graph $AQ_4$ has 12 distinct maximum cliques and that the action of $G:=\Aut(AQ_4)$ on a set of 8 maximum cliques is isomorphic to the action of the automorphism group of two disjoint squares acts on its vertices.  It is then deduced that $\Aut(AQ_4) \cong (D_8 \times D_8) \rtimes C_2$.

\begin{figure}
\centering
\begin{tikzpicture}[scale=2]

\vertex[fill] (ve) at (0,0) [label=below:$e$] {};

\vertex[fill] (ve1) at (2,-1) [label=left:$e_1$] {};
\vertex[fill] (ve2) at (2,0) [label=left:$e_2$] {};
\vertex[fill] (ve3) at (2,1) [label=left:$e_3$] {};

\vertex[fill] (v1111) at (3,-1.5) [label=right:$1111$] {};
\vertex[fill] (v0111) at (3,-0.5) [label=right:$0111$] {};
\vertex[fill] (v0011) at (3,0.5) [label=right:$0011$] {};
\vertex[fill] (ve4) at (3,1.5) [label=right:$e_4$] {};

\node at (3.5,0.5) {$\cdots$};

\path
    (ve1) edge (v1111)
    (ve1) edge (v0111)

    (ve2) edge (v0111)
    (ve2) edge (v0011)

    (ve3) edge (v0011)
    (ve3) edge (ve4)

    (v1111) edge (v0111)
    (v0111) edge (v0011)
    (v0011) edge (ve4)
;
\end{tikzpicture}
\caption{The identity vertex $e$ and its neighbors in $AQ_4$.} 
\label{fig:AQ4:neighborhood:of:e}
\end{figure}
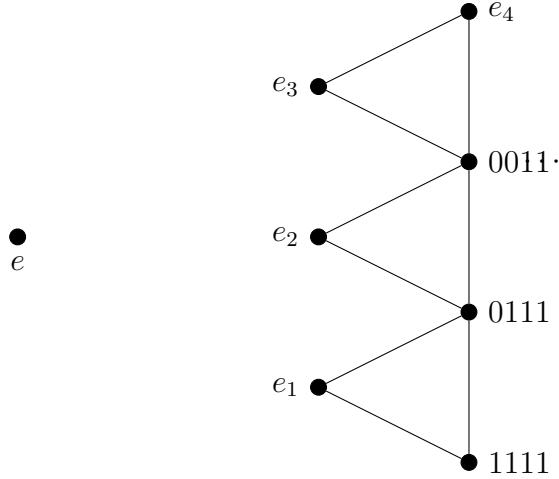

Because $X$ is vertex-transitive, there exists a maximum clique in the graph which contains the identity vertex $e$.  Observe from Figure~\ref{fig:AQ4:neighborhood:of:e} that the size of a maximum clique in the induced subgraph $X[X_1(e)]$ is 3.  Hence the clique number $\omega(X) = 4$.

There are exactly 3 4-cliques in $AQ_4$ containing the identity vertex $e$ (see Figure~\ref{fig:AQ4:neighborhood:of:e}).  Hence, there are a total of $3 \times 16 = 48$ 4-cliques in $AQ_4$; however, in this count, each clique was counted 4 times, and so there are exactly 12 distinct maximum cliques in $AQ_4$, each of size 4.

Let $\mathcal{C}_1 := \{e,e_3,e_4,0011\}$ be the set of 4 vertices of $AQ_4$ that form the maximum clique corresponding to the upper triangle in Figure~\ref{fig:AQ4:neighborhood:of:e}.  We call $\mathcal{C}_1$ an \emph{upper} \emph{clique} of $AQ_4$.  Note that $\mathcal{C}_1$ is also a subspace of $\mathbb{Z}_2^4$ of dimension 2. Let $\mathcal{C}_2 := \mathcal{C}_1 + 1000 = \{e_1,1010,1001,1011\}$ be the coset obtained by translating $\mathcal{C}_1$ by $1000$.  Let $\mathcal{C}_3 := \mathcal{C}_1 + 0100 = \{e_2, 0110, 0101, 0111 \}$ and $\mathcal{C}_4 := \mathcal{C}_1 + 1100 = \{1100, 1110, 1101, 1111\}$.  The set of upper cliques of $AQ_4$ is defined to equal $\{\mathcal{C}_1, \mathcal{C}_2, \mathcal{C}_3, \mathcal{C}_4 \}$.  The set of translations used to obtain these 4 4-cliques is $\{e, 1000, 0100, 1100\}$, a subspace of dimension 2 orthogonal to $\mathcal{C}_1$.  Hence, the set of four upper cliques is a partition of the vertex set $\mathbb{Z}_2^4$ of $AQ_4$. 

Similarly, a lower clique of $AQ_4$ is a 4-clique obtained by a translation of the clique $\mathcal{C}_5 := \{e, e_1, 0111, 1111\}$ corresponding to the lower triangle in Figure~\ref{fig:AQ4:neighborhood:of:e}.  The set of lower cliques is $\{ \mathcal{C}_5, \mathcal{C}_6, \mathcal{C}_7, \mathcal{C}_8 \}$, where $\mathcal{C}_6 := \mathcal{C}_5 + e_4 = \{e_4, 1001, 0110,  1110\}$, $\mathcal{C}_7 := \mathcal{C}_5 + 0011 = \{0011, 1011, e_2, 1100\}$, and $\mathcal{C}_8 := \mathcal{C}_5 + e_3 = \{ e_3, 1010, 0101,  1101\}$. 

The set of middle cliques of $AQ_4$ of course consists of the clique $\mathcal{C}_9 := \{e,0011, e_2, 0111\}$ corresponding to the middle triangle in Figure~\ref{fig:AQ4:neighborhood:of:e} and its translations $\mathcal{C}_{10} := \mathcal{C}_9 + e_4 = \{e_4, e_3, 0101, 0110\}$, $\mathcal{C}_{11} := \mathcal{C}_9 + e_1 = \{e_1, 1011, 1100, 1111 \}$, and $\mathcal{C}_{12} := \mathcal{C}_9 + 1001 = \{1001, 1010, 1101, 1110 \}$.

We have defined above the set $\{ \mathcal{C}_1, \mathcal{C}_2, \ldots,  \mathcal{C}_{12} \}$ of all maximum cliques of $AQ_4$ in terms of translations (cosets) of 2-dimensional subspaces. This yields three different partitions of the vertex set of $AQ_4$, namely $\{ \mathcal{C}_1, \ldots, \mathcal{C}_4\}$, $\{\mathcal{C}_5, \ldots, \mathcal{C}_8 \}$, and $\{ \mathcal{C}_9, \dots, \mathcal{C}_{12} \}$.  Only two of these three partitions are equivalent under the action of the automorphism group of $AQ_4$:
                                                                                                                                                                                                                                                                                                                                                                        
                                                                                                                                                                                                                                                                                                                                                                     \begin{Lemma} \label{lemma:AQ4:two:orbits}                                                                                                                                                                                                                                                                                                                                                                       Let $G:=\Aut(AQ_4)$, and let $\mathcal{C}_1, \ldots, \mathcal{C}_{12}$ be the 12 distinct maximum cliques of $AQ_4$ defined above.  Then, the action of $G$ on $\{\mathcal{C}_1, \ldots, \mathcal{C}_{12} \}$ has exactly 2 orbits: $ \{ \mathcal{C}_1, \ldots, \mathcal{C}_8 \}$ and $\{ \mathcal{C}_9, \ldots, \mathcal{C}_{12} \}$.                                                                                                                                                                                                                                                                                                                                                                       \end{Lemma}

\noindent \emph{Proof}:
Each of the upper cliques $\mathcal{C}_1, \ldots, \mathcal{C}_4$ is a translation of the upper clique $\mathcal{C}_1$.  Because translations are automorphisms of the graph, all 4 upper cliques $\mathcal{C}_1, \ldots, \mathcal{C}_4$ lie in the same orbit of the action of $G$ on maximum cliques.  Similarly, the lower cliques $\mathcal{C}_5, \ldots, \mathcal{C}_8$ lie in the same $G$-orbit, and the middle cliques $\mathcal{C}_9, \ldots, \mathcal{C}_{12} $ lie in the same $G$-orbit.  

The proof of Proposition~\ref{prop:normal:AQn} shows that there exists an element in $G_e$ that takes the upper triangle (of Figure~\ref{fig:AQ4:neighborhood:of:e}) to the lower triangle.  Hence, the upper clique $\mathcal{C}_1$ and the lower clique $\mathcal{C}_5$ lie in the same $G$-orbit, whence $\{C_1, \ldots, C_8\}$ lies in a single $G$-orbit.  

It now suffices to show that the middle clique $\mathcal{C}_9$ and upper clique $\mathcal{C}_1$ do not lie in the same $G$-orbit. The middle clique $\mathcal{C}_9 = \{e, 0011, e_2, 0111\}$ is a 2-dimensional subspace and is closed under addition. Hence, given any vertex $x \in \mathcal{C}_9$, there exists a translation in $\{ \rho_z: z \in \mathcal{C}_9 \}$ that takes $x$ to the identity vertex $e$ and that maps $\mathcal{C}_9$ to itself. Thus, if $\exists g \in G$ such that $g: \mathcal{C}_9 \mapsto \mathcal{C}_1$, then $\exists g \in G_e$ such that $g: \mathcal{C}_9 \mapsto \mathcal{C}_1$. But this is impossible since none of the 8 elements of $G_e \cong D_8$ take the middle triangle to the lower triangle (cf. Figure~\ref{fig:AQ4:neighborhood:of:e}).  Hence, the middle clique $\mathcal{C}_9$ and the upper clique $\mathcal{C}_1$ lie in different $G$-orbits.
\qed

The argument above for $AQ_4$ can be generalized to show that the graph $AQ_n$ $(n \ge 4$) has exactly $(n-1) \times \frac{2^n}{4} = (n-1)2^{n-2}$ distinct maximum cliques, and that the action of $\Aut(AQ_n)$ on this set of maximum cliques has exactly $\lfloor n/2 \rfloor$ orbits.  In particular, it is easy to see that the 4-cliques corresponding to the first triangle and the last triangle of Figure~\ref{fig:induced:subgraph:X1e:AQn} are in the same orbit, the 4-cliques corresponding to the second triangle and the last-but-one triangle are in the same orbit, and so on.

We showed in Lemma~\ref{lemma:AQ4:two:orbits} that the set $\{\mathcal{C}_1, \ldots, \mathcal{C}_8 \}$ of 8 cliques of $AQ_4$ forms a single orbit in the action of $G$ on all maximum cliques.  We study this restricted action of $G$ (i.e. the transitive constituent) further:

\begin{Proposition}  \label{prop:G:action:on:cliques:AQ4:faithful}
 Let $\mathcal{C}_1, \ldots, \mathcal{C}_8$ be the maximum cliques of $AQ_4$ defined above.  Then, the action of $G:=\Aut(AQ_4)$ on $\{\mathcal{C}_1, \ldots, \mathcal{C}_8 \}$ is faithful.
\end{Proposition}

\noindent \emph{Proof}:
Suppose $g \in G$ and $g$ maps each $\mathcal{C}_i$ to itself $(i=1,2,\ldots,8)$.  It suffices to show $g=1$.  Because $g$ maps the lower clique $\mathcal{C}_5$ to itself and the upper clique $\mathcal{C}_1$ to itself, $g$ maps their intersection $\{e\}$ to itself. Hence $g$ fixes the vertex $e$.  Similarly, $g$ fixes $\mathcal{C}_5 \cap \mathcal{C}_2 = \{ e_1 \}$ and $g$ fixes $\mathcal{C}_1 \cap \mathcal{C}_8 = \{e_3\}$.  But it is clear from Figure~\ref{fig:AQ4:neighborhood:of:e} that the only element from $G_e$ that fixes $e_1$ and $e_3$ is the trivial automorphism, i.e. $g=1$.
\qed

If a group $G$ acts on a set $\Omega$, then a block of $G$ is defined to be a subset $\Delta \subseteq \Omega$ such that $\Delta \cap \Delta^g = \Delta$ or $\Delta \cap \Delta^g = \phi$, for all $g \in G$.

\begin{Proposition} \label{prop:AQ4:K:fixed:block}
 Let $\mathcal{K}:=\{\mathcal{C}_1, \mathcal{C}_2, \mathcal{C}_3, \mathcal{C}_4 \}$ and $\mathcal{K}' := \{ \mathcal{C}_5, \mathcal{C}_6, \mathcal{C}_7, \mathcal{C}_8 \}$ denote the set of 4 upper cliques and the set of 4 lower cliques of $AQ_4$, respectively.  Then $\mathcal{K}$ (and hence $\mathcal{K}'$) is a block of the action of $G$ on $\{\mathcal{C}_1, \ldots, \mathcal{C}_8\}$.  
\end{Proposition}

\begin{figure}
\centering
\begin{tikzpicture}[scale=1.1]

\vertex[fill] (v0000) at (2,0) [label=below:$e$] {};
\vertex[fill] (v0001) at (3,0) [label=below:$e_4$] {};
\vertex[fill] (v0010) at (4,0) [label=below:$e_3$] {};
\vertex[fill] (v0011) at (5,0) [label=below:$0011$] {};
\node at (3.5,-1) {$\mathcal{C}_1$};

\vertex[fill] (v1100) at (0,2) [label=left:$1100$] {};
\vertex[fill] (v1101) at (0,3) [label=left:$1101$] {};
\vertex[fill] (v1110) at (0,4) [label=left:$1110$] {};
\vertex[fill] (v1111) at (0,5) [label=left:$1111$] {};
\node at (-2,3.5) {$\mathcal{C}_4$};

\vertex[fill] (v0100) at (2,7) [label=above:$e_2$] {};
\vertex[fill] (v0101) at (3,7) [label=above:$0101$] {};
\vertex[fill] (v0110) at (4,7) [label=above:$0110$] {};
\vertex[fill] (v0111) at (5,7) [label=above:$0111$] {};
\node at (3.5,8) {$\mathcal{C}_3$};

\vertex[fill] (v1000) at (7,2) [label=right:$e_1$] {};
\vertex[fill] (v1001) at (7,3) [label=right:$1001$] {};
\vertex[fill] (v1010) at (7,4) [label=right:$1010$] {};
\vertex[fill] (v1011) at (7,5) [label=right:$1011$] {};
\node at (9,3.5) {$\mathcal{C}_2$};

\path
    (v0000) edge (v1111)
    (v0000) edge (v0100)
    (v0000) edge (v0111)
    (v0000) edge (v1000)
;

\path[dotted]
    (v0001) edge (v1110)
    (v0001) edge (v0101)
    (v0001) edge (v0110)
    (v0001) edge (v1001)
;

\path   
    (v0010) edge (v1101)
    (v0010) edge (v0101)
    (v0010) edge (v0110)
    (v0010) edge (v1010)
;

\path[dotted]
    (v0011) edge (v1100)
    (v0011) edge (v0100)
    (v0011) edge (v0111)
    (v0011) edge (v1011)
;

\end{tikzpicture}
\caption{A partition of the vertex set of $AQ_4$ into cosets of clique $\mathcal{C}_1$.} 
\label{fig:AQ4:4:upper:cliques}
\end{figure}
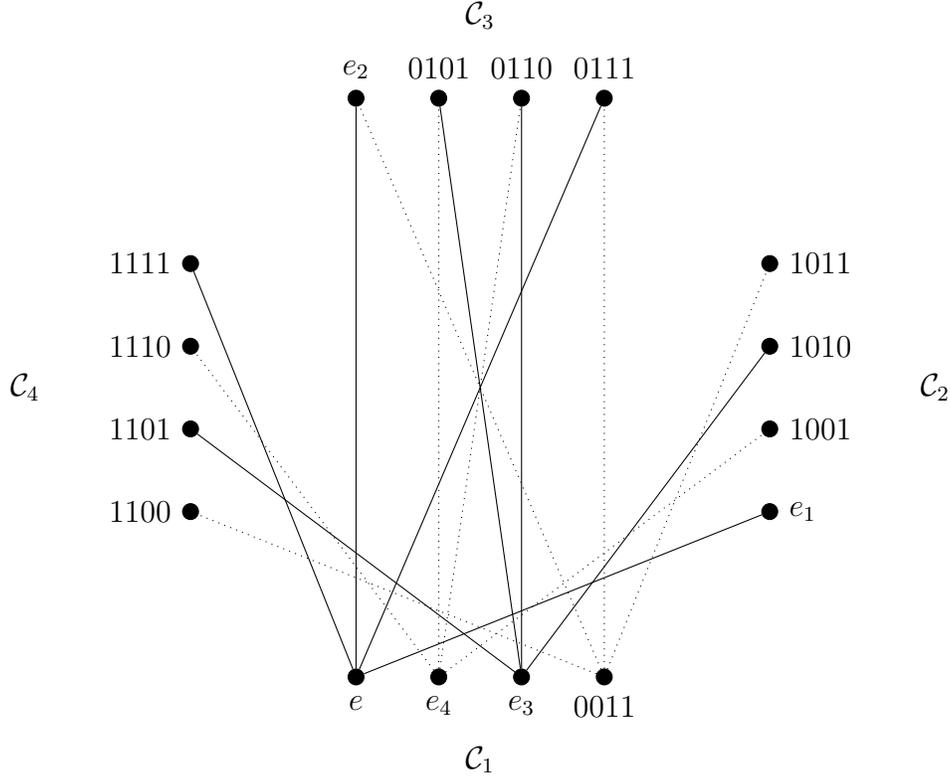

\noindent \emph{Proof}:
Consider the partial drawing of $AQ_4$ shown in Figure~\ref{fig:AQ4:4:upper:cliques}.  In this drawing, the 4 vertices in any upper clique are grouped together and the edges between vertices in the same upper clique are not drawn.  Also, only all edges from the upper clique $\mathcal{C}_1 = \{e,e_3,e_4,0011\}$ to vertices in the other three upper cliques are shown. The edges incident to consecutive vertices of $\mathcal{C}_1$ alternate between solid and dotted format, in order to distinguish between these edges more clearly. Note that the edges between $\mathcal{C}_1$ and $\mathcal{C}_3$ are twice as many as the edges between $\mathcal{C}_1$ and any other upper clique.   Thus, $\{\mathcal{C}_1, \mathcal{C}_3 \}$ is a block of $G$ acting on $\mathcal{K}$.  Studying the automorphism groups of graphs also helps to construct new drawings of graphs.

Observe from Figure~\ref{fig:AQ4:4:upper:cliques} that the neighbors of $e \in \mathcal{C}_1$ in the three other upper cliques are $1111, e_2, 0111$ and $e_1$.  Three of these neighbors, namely $e_1, 0111$ and $1111$, along with $e$, form the lower clique $\mathcal{C}_5 = \{e, e_1, 0111, 1111\}$, whose 4 vertices lie in distinct upper cliques. In other words, there exists a lower clique $\mathcal{C}_5$ that intersects all 4 upper cliques $\mathcal{C}_1, \mathcal{C}_2, \mathcal{C}_3, \mathcal{C}_4$. 

Recall that $G$ acts transitively on $\mathcal{K}$.  Also, an automorphism of the graph must take pairwise disjoint cliques $\mathcal{C}_1, \ldots, \mathcal{C}_4$ to pairwise disjoint cliques.  Hence, if $g \in G$ takes an upper clique, $\mathcal{C}_1$ say, to a lower clique $\mathcal{C}_5$, then $g$ must take each upper clique to a lower clique, i.e. $g$ must interchange $\mathcal{K}$ and $\mathcal{K}'$. Thus, each $g \in G$ either fixes  $\mathcal{K} = \{\mathcal{C}_1, \mathcal{C}_2, \mathcal{C}_3, \mathcal{C}_4 \}$ setwise, or completely moves $\mathcal{K}$ to a set $\mathcal{K}'$ disjoint from $\mathcal{K}$, i.e. $\mathcal{K}$ is a block of the action of $G$ on $\{\mathcal{C}_1, \ldots, \mathcal{C}_8 \}$.
\qed

In the next result, we essentially prove that the action of $G:=\Aut(AQ_4)$ on the cliques $\mathcal{C}_1, \ldots, \mathcal{C}_8$ is isomorphic to the action of the automorphism group of two disjoint squares on the 8 vertices of the two squares.

\begin{Theorem}
 $\Aut(AQ_4) \cong (D_8 \times D_8) \rtimes C_2$.
\end{Theorem}

\noindent \emph{Proof}:
By Lemma~\ref{lemma:AQ4:two:orbits}, $\{\mathcal{C}_1, \ldots, \mathcal{C}_8 \}$ is an orbit of $G$ acting on maximum cliques.  Consider the action of $G$ on $\{ \mathcal{C}_1, \ldots, \mathcal{C}_8 \}$.  By Proposition~\ref{prop:G:action:on:cliques:AQ4:faithful}, this induced action is faithful.  Hence $G$ is determined uniquely by the permutation group image of this induced action. By Proposition~\ref{prop:AQ4:K:fixed:block}, $\mathcal{K} = \{\mathcal{C}_1, \ldots, \mathcal{C}_4 \}$ and $\mathcal{K}' = \{ \mathcal{C}_5, \ldots, \mathcal{C}_8 \}$ are blocks of $G$. 

Since $e_3 \in \mathcal{C}_1$, the translation $\rho_{e_3} \in G$ maps each upper clique $\mathcal{C}_i \in \mathcal{K}$ (which is a coset of the subspace $\mathcal{C}_1$) to itself.  Hence $\rho_{e_3}$ fixes $\mathcal{K}$ pointwise.  On the other hand, $\rho_{e_3}$ interchanges $\mathcal{C}_5$ and $\mathcal{C}_8$ because $e \in \mathcal{C}_5$ and $e_3 \in \mathcal{C}_8$.  Also, $\rho_{e_3}$ interchanges $\mathcal{C}_6$ and $\mathcal{C}_7$ because $e_4 \in \mathcal{C}_6$ and $e_4+e_3 \in \mathcal{C}_7$. Thus, $\rho_{e_3}$ effects the permutation $(\mathcal{C}_5, \mathcal{C}_8)(\mathcal{C}_6, \mathcal{C}_7) \in \Sym( \{ \mathcal{C}_1, \ldots, \mathcal{C}_8 \} )$.  Similarly, $\rho_{e_4}$ fixes $\mathcal{K}$ pointwise and effects the permutation $(\mathcal{C}_5, \mathcal{C}_6)(\mathcal{C}_7, \mathcal{C}_8)$ on $\mathcal{K}'$. 

Let $A$ be the linear extension of the map $(e_1,e_2,e_3,e_4) \mapsto (e_1,e_2,e_4,e_3)$, i.e. $A$ interchanges the last 2 coordinates of a vector. The coset $\mathcal{C}_1$ is closed under action by $A$ and each vector in the span of $\{e_1,e_2\}$ is also fixed by $A$.  So $A$ maps the coset $\mathcal{C}_1 +a$ $(a \in \mbox{Span}\{e_1,e_2\})$ to itself because $(\mathcal{C}_1+a)^A = \mathcal{C}_1^A + a^A = \mathcal{C}_1 + a$. Thus, $A$ fixes $\mathcal{K}$ pointwise. 

The linear transformation $A$ just defined is invertible.  Note also that $A$ fixes $S$ setwise.  Hence $A \in \Aut(\mathbb{Z}_2^4,S) \subseteq G_e$.  Thus, $A$ permutes the cliques $\mathcal{C}_1, \ldots, \mathcal{C}_8$ among themselves.  Since $A$ fixes $\mathcal{K}$ pointwise, $A$ fixes $\mathcal{K}'$ setwise. It can be verified that $A$ fixes $\mathcal{C}_5$ since $e \in \mathcal{C}_5$.  Also, $A$ interchanges $\mathcal{C}_6$ and $\mathcal{C}_8$ because $e_4 \in \mathcal{C}_6$ and $e_4^A = e_3 \in \mathcal{C}_8$.  Thus, $A$ effects the permutation $(\mathcal{C}_6, \mathcal{C}_8) \in \Sym(\{\mathcal{C}_1, \ldots, \mathcal{C}_8 \})$.  

\begin{figure}
\centering
\begin{tikzpicture}[scale=1.5]

\vertex[fill] (v1) at (-1,1) [label=left:$\mathcal{C}_1$] {};
\vertex[fill] (v2) at (1,1) [label=right:$\mathcal{C}_2$] {};
\vertex[fill] (v3) at (1,-1) [label=right:$\mathcal{C}_3$] {};
\vertex[fill] (v4) at (-1,-1) [label=left:$\mathcal{C}_4$] {};
\draw[thick] (v1) -- (v2) -- (v3) -- (v4) -- (v1);
\draw[dashed] (-1.5,1.5) -- (1.5, -1.5);
\node at (1.6, -1.6) {$A'$};
\draw[dashed] (0,1.5) -- (0, -1.5);
\node at (0, -1.6) {$\rho_{e_1}$};
\draw[dashed] (-1.5,0) -- (1.5, 0);
\node at (1.9, 0) {$\rho_{1111}$};
\node at (0, -2.2) {$\mathcal{K}$};

\vertex[fill] (v5) at (4,1) [label=left:$\mathcal{C}_5$] {};
\vertex[fill] (v6) at (6,1) [label=right:$\mathcal{C}_6$] {};
\vertex[fill] (v7) at (6,-1) [label=right:$\mathcal{C}_7$] {};
\vertex[fill] (v8) at (4,-1) [label=left:$\mathcal{C}_8$] {};
\draw[thick] (v5) -- (v6) -- (v7) -- (v8) -- (v5);
\draw[dashed] (3.5,1.5) -- (6.5, -1.5);
\node at (6.6, -1.6) {$A$};
\draw[dashed] (5,1.5) -- (5, -1.5);
\node at (5, -1.6) {$\rho_{e_4}$};
\draw[dashed] (3.5,0) -- (6.5, 0);
\node at (6.9, 0) {$\rho_{e_3}$};
\node at (5, -2.2) {$\mathcal{K}'$};

\end{tikzpicture}
\caption{$G:=\Aut(AQ_4)$ acting on maximum cliques.}
\label{fig:AQ4:cliques:square}
\end{figure}
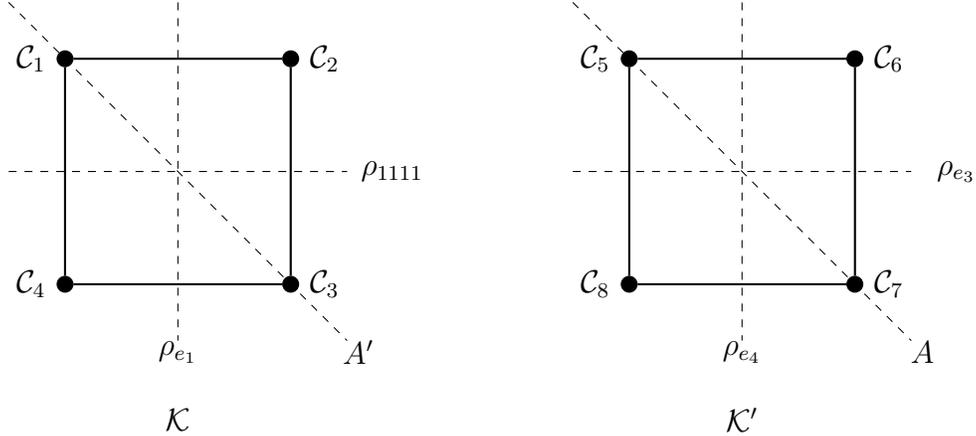

Thus $N_1 := \langle \rho_{e_3}, \rho_{e_4}, A \rangle$ $=\langle (\mathcal{C}_5, \mathcal{C}_8)(\mathcal{C}_6, \mathcal{C}_7), (\mathcal{C}_5, \mathcal{C}_6)(\mathcal{C}_7, \mathcal{C}_8), (\mathcal{C}_6, \mathcal{C}_8) \rangle \cong D_8$ (see Figure~\ref{fig:AQ4:cliques:square}).

Similarly, consider $N_2:=\langle \rho_{e_1}, \rho_{1111}, A' \rangle$, where $\rho_{e_1}$ maps each lower clique to itself and effects the permutation $(\mathcal{C}_1, \mathcal{C}_2)(\mathcal{C}_3,\mathcal{C}_4)$ on the upper cliques; $\rho_{1111}$ effects the permutation $(\mathcal{C}_1, \mathcal{C}_4)(\mathcal{C}_2, \mathcal{C}_3)$.  Let $A'$
be the linear transformation that takes $e_1$ to $1111$ and $e_i$ to itself $(i=2,3,4)$.  Then, $A'$ fixes $S$ setwise, whence $A' \in \Aut(\mathbb{Z}_2^4,S) \subseteq G_e$. It can be verified that $A'$ effects the permutation $(\mathcal{C}_2, \mathcal{C}_4) \in \Sym(\{\mathcal{C}_1, \ldots, \mathcal{C}_8 \})$.  Thus, $N_2 \cong D_8$.

The permutations in $N_1$ fix each element of $\mathcal{K}$ and move only elements in $\mathcal{K}'$.  The permutations in $N_2$ fix each element of $\mathcal{K}'$ and move only elements in $\mathcal{K}$.  Hence, the elements of $N_1$ and $N_2$ commute with each other and $N_1 \cap N_2=1$. Thus, $N_1 N_2 = N_1 \times N_2 \cong D_8 \times D_8$. 

A permutation of $\mathcal{C}_1, \ldots, \mathcal{C}_8$ that interchanges $\mathcal{K}$ and $\mathcal{K}'$ is one induced by the linear map $B: (e_1,e_2,e_3,e_4) \mapsto (e_3, e_2, e_1, 1111)$ (cf. Figure~\ref{fig:AQ4:neighborhood:of:e}). 
Then $B$ interchanges the upper and lower cliques.  The subgroup $\langle B \rangle \cong C_2$ normalizes $N_1 \times N_2$, which has index 2 in $(N_1 \times N_2) \langle B \rangle $.  Thus, $(D_8 \times D_8) \rtimes C_2$ is isomorphic to a subgroup of $G$. By order considerations (cf. Theorem~\ref{thm:autgroup:sdp:AQn}), $(D_8 \times D_8) \rtimes C_2$ is the full automorphism group $G$.
\qed

%--------------------------

%---------------------------------------------------

\section{Nontrivial blocks of $AQ_4$}

In this section we obtain an equivalent condition for a subset of vertices of $AQ_n$ to be a nontrivial block of the graph, and we determine all the nontrivial blocks of $AQ_4$.  We first recall some definitions on blocks (cf. Wielandt \cite{Wielandt:1964}, Biggs \cite{Biggs:1993}).

Let $G \le \Sym(V)$ be a transitive permutation group.  A block of $G$ is a subset $\Delta \subseteq V$ that satisfies the condition $\Delta \cap \Delta^g = \Delta$ or $\Delta \cap \Delta^g = \phi$, for all $g \in G$.  
A block $\Delta$ of $G$ is said to be nontrivial if $1 < |\Delta| < |V|$.  The set $\{\Delta^g: g \in G\}$ of all translates of a block is called a complete block system (or a system of imprimitivity) and forms a partition of $V$.   A permutation group $G$ is said to be imprimitive if it has a nontrivial block. 

A graph $X=(V,E)$ is imprimitive if it is transitive and its automorphism group is imprimitive.  The blocks of $X$ are defined to be the blocks of its automorphism group. Thus, if $W \subseteq V$ is a block of $X$ and $g \in \Aut(X)$, then $W^g$ and $W$ are either equal or disjoint (i.e. there cannot be overlap that is only partial).  For example, for the hypercube graph, an antipodal pair $\{x, \overline{x} \}$ of vertices is a nontrivial block. The hypercube is bipartite, and the set of all even weight vertices is another nontrivial block. It can be shown (cf. Biggs \cite[Proposition 22.3]{Biggs:1993}) that the hypercube has only two nontrivial block systems. 

We will use the following result: 

\begin{Theorem} \cite[Theorem 7.4, 7.5]{Wielandt:1964} \label{thm:lattices:blocks:subgroups}
 Let $G \le \Sym(V)$ be a transitive permutation group and let $e \in V$.  Then, the lattice of blocks of $G$ which contain $e$ is isomorphic to the lattice of subgroups of $G$ that contain $G_e$. More specifically, the correspondence between blocks $\Delta$ and subgroups $H$ is as follows: $\Delta := e^H$ (the orbit of $e$ under the action of $H$), and $H:=G_{ \{ \Delta \}}$ (the set of elements of $G$ that fixes $\Delta$ setwise).  
\end{Theorem}

To obtain the blocks of the augmented cube graph $AQ_n$, we determine the subgroups of $G:=\Aut(AQ_n)$ that contain $G_e$ and then use Theorem~\ref{thm:lattices:blocks:subgroups}.   For $z \in \mathbb{Z}_2^n$, let $(\rho_z: \mathbb{Z}_2^n \rightarrow \mathbb{Z}_2^n, x \mapsto x+z)$ denote translation by $z$. For a subset $K \subseteq \mathbb{Z}_2^n$, let $R(K):= \{ \rho_z: z \in K\}$.  We say that a subset $K \subseteq \mathbb{Z}_2^n$ is closed under action by $G_e$ if $k^g \in K$, for all $k \in K, g \in G_e$. 

\begin{Lemma} \label{lemma:condition:on:subset:K}
 Let $X:=\Cay(\mathbb{Z}_2^n,S)$ be a normal Cayley graph (here, $X$ need not be $AQ_n$), and let $G:=\Aut(X) = R(\mathbb{Z}_2^n) \rtimes G_e$.  Let $K$ be a subset of vertices of the graph which forms a subspace of $\mathbb{Z}_2^n$ and which is closed under action by $G_e$.  Then $R(K)G_e$ is a subgroup of $G$.
\end{Lemma}

\noindent \emph{Proof}:
The group $G$ is finite, so it suffices to show that $R(K)G_e$ satisfies closure.  We need to show that for all $z_1,z_2 \in K$ and $g_1,g_2 \in G_e$, there exist $z_3 \in K$ and $g_3 \in G_e$ such that $(\rho_{z_1} g_1)(\rho_{z_2} g_2) = \rho_{z_3} g_3$.  This latter equation is equivalent to the condition $x^{(\rho_{z_1} g_1)(\rho_{z_2} g_2)} = x^{\rho_{z_3} g_3}$, for all $x \in V(X)$.  The left hand side is $x^{{\rho_{z_1} g_1}{\rho_{z_2} g_2}} = (x+z_1)^{g_1 \rho_{z_2} g_2}$ $=(x^{g_1} + z_1^{g_1} + z_2)^{g_2}$ $=x^{g_1 g_2} + z_1^{g_1 g_2} + z_2^{g_2}$.  
Here, we used the fact that $g_i$ is a linear transformation, which is the case because for a normal Cayley graph $G_e \subseteq \Aut(\mathbb{Z}_2^n)$. The right hand side is $x^{\rho_{z_3} g_3} = x^{g_3} + z_3^{g_3}$.  So we get that
\[
 x^{g_1 g_2} + z_1^{g_1 g_2} + z_2^{g_2} = x^{g_3}+z_3^{g_3}, \forall x \in V(X).
\]
Take $g_3=g_1 g_2$. Then $z_3$ must satisfy the equation $z_1^{g_1 g_2} + z_2^{g_2} = z_3^{g_1 g_2}$, or equivalently, the equation $z_1 + z_2^{g_1^{-1}} = z_3$.  By hypothesis, $K$ is closed under addition and under action by $g_1 \in G_e$.  Hence, there exists a $z_3 \in K$ satisfying the equation. It follows that $R(K)G_e$ is a subgroup of $G$.
\qed

We now show that the subgroups $R(K)G_e$ obtained from Lemma~\ref{lemma:condition:on:subset:K} are all the subgroups of $G$ that contain $G_e$:

\begin{Proposition}
 Let $X:=\Cay(\mathbb{Z}_2^n,S)$ be a Cayley graph and let $G:=\Aut(X)$.  If $G_e \le H \le G$, then $H=R(K) G_e$ for some subspace $K$ of $\mathbb{Z}_2^n$ that is closed under action by $G_e$.
\end{Proposition}

\noindent \emph{Proof}:
Suppose $H$ satisfies $G_e \le H \le G$. Then $H$ can be expressed as a disjoint union of left cosets of $G_e$ in $H$, i.e.  $H = h_1 G_e \cup h_2 G_e \cup \ldots \cup h_k G_e$, where $k$ is the index of $G_e$ in $H$.  Note that $h_i \in H \le G = R(\mathbb{Z}_2^n) G_e$, and so each $h_i$ can be expressed as $h_i = \rho_{z_i} g_i$ for some $z_i \in \mathbb{Z}_2^n, g_i \in G_e$. Then, $h_i G_e = \rho_{z_i} g_i G_e = \rho_{z_i} G_e$. Thus, $H=R(K)G_e$, where $K:=\{z_1,\ldots,z_k\}$.  This proves that if $H$ is a subgroup of $G$ which contains $G_e$, then $H=R(K) G_e$ for some subset $K \subseteq \mathbb{Z}_2^n$.

To show $K$ is a subspace, suppose $z_1,z_2 \in K$ (it is clear $K$ is nonempty). We show $z_1+z_2 \in K$.  By closure in the subgroup $R(K)G_e$,  $(\rho_{z_1} g_1)(\rho_{z_2} g_2) \in R(K) G_e$, for all $g_1, g_2 \in G_e$. In particular, taking $g_1 = g_2 = 1$, we get that $\rho_{z_1} \rho_{z_2} \in R(K) G_e$. This implies that there exist $z_3 \in K$ and $g_3 \in G_e$ such that $\rho_{z_1 + z_2} = \rho_{z_3} g_3$. The right hand side takes vertex $z_3$ to $z_3^{\rho_{z_3} g} = e^g =e$, whence $z_3^{\rho_{z_1+z_2}} = e$. But this implies that $z_3 = z_1 + z_2$.  Thus, for closure to hold in $R(K)G_e$, $R(K)$ must also contain $\rho_{z_1+z_2}$, i.e. $K$ is closed under addition and hence is a subspace.

It remains to show that $K$ is closed under action by $G_e$. Let $z_1 \in K, g_1 \in G_e$.  We show $z_1^{g_1} \in K$.  Since $R(K)G_e$ is a subgroup, there exist $z_2 \in K$ and $g_2 \in G_e$ such that $\rho_{z_1} g_1 = g_2 \rho_{z_2}$, i.e. $x^{\rho_{z_1} g_1 \rho_{z_2}} = x^{g_2}$, for all $x \in V(X)$.  Taking $x=e$, we get $z_1^{g_1} = z_2$, i.e. $z_1^{g_1}$ is also in $K$. 
\qed

 Thus, the subgroups $H$ of $G$ that contain $G_e$ are in bijective correspondence with the subspaces $K$ of $\mathbb{Z}_2^n$ that are closed under action by $G_e$.  The conditions of Lemma~\ref{lemma:condition:on:subset:K} require $K$ to be closed under action by $G_e$. This implies $K$ is a union of orbits of $G_e$.  

Suppose we have a subset $K \subseteq \mathbb{Z}_2^n$ satisfying the conditions of Lemma~\ref{lemma:condition:on:subset:K}. By Theorem~\ref{thm:lattices:blocks:subgroups}, the subgroup $H:=R(K)G_e$ gives rise to the block $\Delta = e^H = e^{R(K)G_e} = K^{G_e} = K$.  Given such a block $\Delta$, a complete block system $\{\Delta^g: g \in G\}$ is the set  of all translations of $\Delta$: if $\rho_z g \in G$ $(z \in \mathbb{Z}_2^n, g \in G_e)$, then $\Delta^{\rho_z g} = K^{\rho_z g} = (K+z)^g = K^g + z^g = K+z^g = \Delta+z^g$. Thus, to determine all nontrivial block systems of a normal Cayley graph $X=\Cay(\mathbb{Z}_2^n,S)$, it suffices to determine all subspaces $K$ of $\mathbb{Z}_2^n$ that are closed under action by $G_e$. 

By the arguments above, we get the following characterization for the nontrivial blocks of a normal Cayley graph:

\begin{Theorem} \label{thm:characterization:blocks}
 Let $X:=\Cay(\mathbb{Z}_2^n,S)$ be a normal Cayley graph and let $G:=\Aut(X)$.  Then, the nontrivial blocks of $X$ that contain the identity vertex $e$ are precisely the subspaces of $\mathbb{Z}_2^n$ that are closed under action by $G_e$.
\end{Theorem}

We now apply the above results to the special case where the graph is $AQ_4$.

\begin{Proposition} \label{prop:nontrivial:blocks:AQ4}
Each of the following subsets is a nontrivial block of $AQ_4$:
$\Delta = \{e,e_2\}$; $\Delta' := \{e, e_2, 0011, 0111\}$; $\Delta'' := \{ e, e_2, 0011, 0111,$ $1001, 1010, 1101, 1110 \}$. 
\end{Proposition}

\noindent \emph{Proof}:
The vertex $e_2$ is fixed by $G_e$ (cf.  Figure~\ref{fig:AQ4:neighborhood:of:e}).  Hence, the subspace spanned by $\{e_2\}$ is fixed by $G_e$. By Theorem~\ref{thm:characterization:blocks}, this subspace is a block of $AQ_4$.  The four vectors in $\Delta'$ form a 2-dimension subspace $\mbox{Span} \{e_2,0011\}$ of $\mathbb{Z}_2^4$ and this subset is closed under action by $G_e$ (see Figure~\ref{fig:AQ4:neighborhood:of:e}).  By Theorem~\ref{thm:characterization:blocks}, $\Delta'$ is a block of $AQ_4$.

It suffices to show that $\Delta''$ is a union of orbits of $G_e$ and that $\Delta''$ forms a subspace.  It is clear from the drawing of $AQ_4$ (see Figure~\ref{fig:AQ4:neighborhood:of:e}) that $\{e\}$, $\{e_2\}$ and $\{0011, 0111\}$ are fixed blocks (hence union of orbits) of $G_e$. 
Since $G_e$ fixes $e_2$, $G_e$ fixes the set $\{1001, 1010, 1101, 1110\}$ of non-neighbors of $e_2$ in $X_2(e)$ setwise. Thus, this set is a fixed block of $G_e$ and $\Delta''$ is a union of orbits of $G_e$. It can be verified that $\Delta''$ is a subspace of dimension 3 spanned by $\{0011, 0100, 1001\}$.  Hence $\Delta'' \subseteq \mathbb{Z}_2^4$ is such that $R(\Delta'') G_e$ is a subgroup, and hence, the corresponding set $\Delta''$ is a block of $G$.
\qed

\begin{Theorem}
 The 3 nontrivial blocks given in Proposition~\ref{prop:nontrivial:blocks:AQ4} are the only nontrivial blocks of $AQ_4$ which contain the vertex $e$.
\end{Theorem}

\noindent \emph{Proof}:
By Theorem~\ref{thm:characterization:blocks}, it suffices to find all the subspaces $K$ of $\mathbb{Z}_2^4$ that are closed under action by $G_e$.  We can find the orbits of $G_e$, and determine which fixed blocks (union of orbits) of $G_e$ are subspaces.  Since $K$ must be nontrivial, $K \ne \{e\}$ and $K \ne \mathbb{Z}_2^4$. 

It is clear from Figure~\ref{fig:AQ4:neighborhood:of:e} that $\Delta_1 := \{e\}$, $\Delta_2 := \{e_2\}$, $\Delta_3 := \{0011, 0111\}$ and $\Delta_4 := \{e_1, 1111, e_3, e_4 \}$ are orbits of $G_e$. Note that $\Delta_4$ spans $\mathbb{Z}_2^4$, and so $K$ does not contain $\Delta_4$.  Define the linear transformations $A_1: (e_1,e_2,e_3,e_4): \mapsto (e_1,e_2,e_4,e_3)$, $A_2: (e_1,e_2,e_3,e_4) \mapsto (1111, e_2, e_3, e_4)$, and $A_3: (e_1,e_2,e_3,e_4) \mapsto (e_3,e_2,e_1,1111)$.  Then $\langle A_1,A_2,A_3 \rangle = G_e \cong D_8$. 

Let $x:=0101$. Then $x^{A_1} = 0110$, $x^{A_3} = 1011$ and $x^{A_1 A_3} = 1100$.  Hence, $\Delta_5 := \{0101, 0110, 1011, 1100 \}$ lies in a single $G_e$-orbit.  Let $y:=1110$. Then $y^{A_1} = 1101$, $y^{A_2} = 1001$ and $y^{A_1 A_2} = 1010$.  Hence, $\Delta_6 := \{1110, 1101, 1001, 1010 \}$ lies in a single $G_e$-orbit. If $\Delta_5 \dot{\cup} \Delta_6 \subseteq K$, then $|K| \ge |\Delta_5| + |\Delta_6| + |\{e\}| = 9$, but since $K$ is a  subspace, $|K|$ is a power of 2, and so $K = \mathbb{Z}_2^4$, a contradiction. Thus, at most one of $\Delta_5$ or $\Delta_6$ is contained in $K$.  If $K \supseteq \Delta_5$, then $|K| \ge 5$, and so $|K|=8$. But it can be verified that $\Delta_1 \dot{\cup} \Delta_2 \dot{\cup} \Delta_3 \dot{\cup} \Delta_5$ is not a subspace.  The subset $\Delta_1 \dot{\cup} \Delta_2 \dot{\cup} \Delta_3 \dot{\cup} \Delta_6$ is the subspace $\Delta''$ of Proposition~\ref{prop:nontrivial:blocks:AQ4}.  If $K$ contains neither $\Delta_5$ nor $\Delta_6$, then the only choices for $K$ are $\Delta_1 \dot{\cup} \Delta_2$ and $\Delta_1 \dot{\cup} \Delta_2 \dot{\cup} \Delta_3$ (denoted $\Delta$ and $\Delta'$, respectively, in Proposition~\ref{prop:nontrivial:blocks:AQ4}).  
\qed

\section{Acknowledgements}

Thanks are due to the authors of \cite{Choudum:Sunitha:2002} \cite{Choudum:Sunitha:2001} for sending copies of their papers.
%---------------------------
{
\bibliographystyle{plain}
\bibliography{refsaut}
}
\end{document}